\documentclass[11pt]{amsart}
\usepackage{amssymb,amsthm}
\usepackage{mathtools}
\usepackage{graphicx}
\usepackage[left=1.125in,right=1.125in,top=1.25in,bottom=1.25in]{geometry}
\usepackage{verbatim}
\usepackage{hyperref}
\usepackage{charter,eulervm}
\usepackage{csquotes}
\usepackage{color}
\usepackage{tikz}
\usetikzlibrary{cd}
\usepackage{setspace}
\usepackage{todonotes}

\theoremstyle{plain}
\newtheorem{theorem}{Theorem}[section]
\newtheorem*{theorem*}{Theorem}
\newtheorem{corollary}[theorem]{Corollary}
\newtheorem{proposition}[theorem]{Proposition}

\newtheorem{lemma}[theorem]{Lemma}

\theoremstyle{definition}
\newtheorem{definition}[theorem]{Definition}
\newtheorem*{definition*}{Definition}

\newtheorem{example}[theorem]{Example}

\theoremstyle{remark}
\newtheorem*{remark*}{Remark}

\newtheorem*{notation*}{Notation}

\DeclareMathOperator{\coz}{coz}

\DeclareMathOperator{\iruc}{iruc}
\DeclareMathOperator{\oruc}{oruc}

\newcommand{\downset}[1]{\left\downarrow{#1}\right\downarrow}
\newcommand{\upset}[1]{\left\uparrow{#1}\right\uparrow}

\newcommand\pc{\xrightarrow{\text{\raisebox{-2pt}[-6pt][-6pt]{$\bullet$}}}}
\newcommand{\ol}[1]{\overline{#1}}

\newcommand{\mbf}[1]{\mathbf{#1}}
\newcommand{\mbb}[1]{\mathbb{#1}}
\newcommand{\mcal}[1]{\mathcal{#1}}

\newcommand{\sbw}[2]{\smashoperator[#1]{\bigwedge_{#2}}}
\newcommand{\sbv}[2]{\smashoperator[#1]{\bigvee_{#2}}}

\newcommand{\setof}[2]{\left\{\,#1 : #2\,\right\}}
\newcommand{\ssetof}[2]{\left\{#1\right\}_{#2}}
\newcommand{\qtq}[1]{\quad\text{#1}\quad}
\newcommand{\abs}[1]{\left|#1\right|}
\newcommand{\map}[3]{#1 \colon #2 \to #3}
\newcommand{\exR}{\ol{\mbb{R}}}
\newcommand{\xra}[1]{\xrightarrow{#1}}

\begin{document}
\title[Four uniform completions of an archimedean vector lattice]{The four uniform completions of a \\unital archimedean vector lattice}
\author{Richard N. Ball}
\address[Ball]{Department of Mathematics, University of Denver, Denver, CO 80208, U.S.A.}
\author{Anthony W. Hager}
\address[Hager]{Department of Mathematics, Wesleyan University, Middletown, CT 06459, U.S.A.}
\date{\today}
\thanks{File name: Pointfree relative uniform convergence.tex}
\keywords{pointwise convergence, relative uniform convergence, unital archimedean vector lattice, completely regular frame, $P$-frame}
\subjclass[2020]{06D22,06F20, 18F70, 46A40}

\begin{abstract}
	In the category \(\mbf{V}\) of unital archimedean vector lattices, four notions of uniform completeness obtain. In all cases completeness requires the convergence of uniformly Cauchy sequences; the completions are distinguished by the manner in which the convergence is regulated.
	\begin{itemize}
		\item 
		Ordinary uniform convergence is regulated by the canonical unit \(1\).
		
		\item 
		Inner relative uniform convergence, here termed \emph{iru-convergence}, is regulated by an arbitrary positive element.
		
		\item 
		Outer relative uniform convergence, here termed \emph{oru-convergence}, is regulated by an arbitrary positive element of a vector lattice containing the given object as a sub-vector lattice.
		
		\item 
		\emph{\(*\)-convergence} is equivalent to ordinary uniform convergence on certain specified quotients of the vector lattice. 
	\end{itemize}   
	In each case the complete objects form a full monoreflective subcategory of \(\mbf{V}\), denoted respectively \(\mbf{ucV}\), \(\mbf{irucV}\),  \(\mbf{orucV}\), and \(\mbf{*cV}\). 

	In this article we provide a unified development of these completions by means of a novel pointfree variant of the classical Yosida adjunction. In this adjunction, the geometric counterpart of a given vector lattice \(G\) is a compactification, i.e., a dense frame surjection \(\map{q}{K}{M}\) with compact domain. \(G\) is then represented as a vector lattice of compacification morphisms \(\map{(g', g)}{p}{q}\), where \(\map{p}{\mcal{O}\ol{\mbb{R}}}{\mcal{O}\mbb{R}}\) is the compacification which is the frame counterpart of the inclusion of the reals \(\mbb{R}\) in the extended reals \(\ol{\mbb{R}}\). We require that a compactification morphism make the diagram commute.
	\[
		\begin{tikzcd}
			\mcal{O}\ol{\mbb{R}}		\arrow{r}{g'}
										\arrow{d}[swap]{p}
			&K							\arrow{d}{q}\\
			\mcal{O}\mbb{R}				\arrow{r}[swap]{g}
			&M
		\end{tikzcd}
	\]
	This representation is well suited to our purposes, for it provides insight into the distinctions among the uniform completions of a vector lattice \(G\).  Passage to the ordinary uniform completion of \(G\) does not change the frames \(K\) or \(M\), but does fill in all possibilities for entries \(g'\) in the top arrow of the diagram. Passage to the inner relative uniform completion of \(G\) typically enlarges the frame \(K\) but does not change \(M\), and fills in all possibilities for bounded entries \(g\) in the bottom arrow of the diagram. Passage to the \(*\)-completion does not change \(M\) but enlarges \(K\) to the compact (\v{C}ech-Stone) coreflection of \(M\) and fills in all possibilities for \(g\) and \(g'\). The passage to the outer relative uniform completion of \(G\) is an enormous enlargement. \(M\) is enlarged to its \(P\)-frame reflection and \(K\) is replaced by the compact (\v{C}ech-Stone) coreflection of \(M\), and all possibilities for \(g\) and \(g'\) are included. In fact, this completion is the largest monoreflection of \(G\), namely its epicompletion.
\end{abstract}

\maketitle


\section{Introduction}

Relative uniform convergence must play a central role in any consideration of archimedean vector lattices, for it is precisely the closure of a convex sub-vector lattice with respect to this convergence which characterizes the archimedeanness of the quotient. Likewise the corresponding completions occupy a central position in the theory. There are four such completions of a given archimedean vector lattice. The first of these, classical uniform completion with regulator \(1\), is sufficiently familiar as to receive scant further attention in the article in hand; the Stone Weierstrass Theorem is mentioned only in the proofs of Proposition \ref{Prop:25} and Theorem \ref{Thm:1}. Nor shall we concern ourselves much with the fourth completion, the so called \emph{\(*\)-completion.} Introduced in \cite{BallHager:2006} with a greatly simplified proof of its central result provided in \cite{BallHager:2019}, this completion of a vector lattice \(G\) yields \(\mcal{R}M\), the vector lattice of frame homomorphisms \(\mcal{O}\mbb{R} \to M\), where \(\mcal{O}\mbb{R}\) is the topology of the real numbers \(\mbb{R}\) and \(M\) is the Madden frame of \(G\). This article mainly concerns the second and third completions, here termed respectively the \emph{inner and outer relative uniform completions,} or iru-completion and oru-completion for short. 

Section \ref{Sec:YosThm} outlines the mathematical apparatus which underlies the analysis; it is a general form of the (pointfree) Yosida adjunction (Theorem \ref{Thm:3}) introduced in \cite{Ball:2021}.  A virtue of its format is that it associates with each vector lattice two frames linked by a compactification map \(\map{q}{K}{M}\). Roughly speaking, \(K\) is determined by and diagnostic of the bounded elements of the vector lattice, and \(M\) is determined by and diagnostic of the unbounded elements of the vector lattice. The analysis is necessarily pointfree, for although \(K\) is always spatial, \(M\) need not be.

Given that fact, it is no small irony that pointwise convergence plays an important role in our analysis (Section \ref{Sec:PtwiseCon}). This convergence, which enjoys a graceful articulation in a pointfree setting, is implied by relative uniform convergence (Proposition \ref{Prop:24}). Its importance derives from the fact that it is precisely the pointwise joins which are preserved by all vector lattice homomorphisms (Proposition \ref{Prop:3}(2)). Furthermore, the pointwise Cauchy completion of a vector lattice \(G\) is its epicompletion (\cite{Ball:2018}), which coincides with the oru-completion of \(G\). The epicompletion of \(G\) also coincides with \(\mcal{R}P\), where \(M \to P\) is the \(P\)-frame reflection of \(M\) (\cite{BallWWZenk:2010}), \(M\) is the Madden frame of \(G\), \(\mcal{R}P\) designates the vector lattice of frame homomorphisms \(\mcal{O}\mbb{R} \to P\), and \(\mcal{O}\mbb{R}\) designates the topology of the real numbers. We briefly review the construction of this reflection in the preliminary Section \ref{Sec:Prelim}; our reason for doing so is to emphasize that the key step in the construction, Lemma \ref{Lem:5}, is a crucial ingredient of our analysis of the oru-completion in Proposition \ref{Prop:4}.

In \cite{Veksler:1969} Veksler constructs the iru-completion of a vector lattice inside its Dedekind MacNeille completion in its maximal essential extension. In Theorem \ref{Thm:13} we offer an alternative construction which is not only conceptually simpler, but also makes it clear that the completion is much closer to \(G\) than is implied by Veksler's construction: the iru-completion of \(G\) is just its relative uniform closure in \(\mcal{R}M\). This is an immediate consequence of the fact that \(\mcal{R}M\) is iru-complete (Theorem \ref{Thm:12}).  

Our notation is standard for the most part, save for our use of 
\[
	\downset{A}
	\equiv \setof{b}{\exists a \in A\ (b \leq a)}
\]
to designate the downset generated by a subset $A \subseteq P$, and dually for our use of $\upset{A}$ for the upset generated by $A$. The letters \(G\) and \(H\) signify archimedean vector lattices with distinguished weak order unit designated by \(1\), while the letters \(K\), \(L\), \(M\), and \(N\) signify completely regular frames.

\section{Preliminaries}\label{Sec:Prelim}

In this section we briefly recall the epicompletion of a vector lattice as constructed from the \(P\)-frame reflection of its Madden frame. The construction rests on the notion of freely complementing a frame element in an extension. This idea will surface in the proof of Proposition \ref{Prop:4}. 

\subsection*{Freely complementing frame elements}

\begin{definition*}[\(\map{o_a}{L}{O_a}\), \(\map{c_a}{L}{C_a}\)]
	Let \(a\) be an element of a frame \(L\). The \emph{open quotient associated with \(a\)} is the map 
	\[
		\map{o_a}{L}{\setof{a \to b}{b \in L}}
		\equiv O_a
		=(b \mapsto a \to b).
	\]
	The \emph{closed quotient associated with \(a\)} is the map
	\[
		\map{c_a}{L}{\upset{a}}
		\equiv C_a
		= (b \mapsto a \vee b).
	\]	
\end{definition*}

\begin{definition}[element \(a\) freely complemented in an extension \(\map{k}{L}{K}\)]
	We say that an element \(a\) of a frame \(L\) is \emph{complemented in an extension \(\map{k}{L}{K}\)} if \(k\) is an injective frame homomorphism such that \(K\) contains an element \(b\) complementary to \(k(a)\), i.e., \(k(a) \wedge b = \bot\) and \(k(a) \vee b = \top\). We say that \emph{\(a\) is freely complemented in the extension} if any extension of \(L\) in which \(a\) is complemented factors uniquely through \(k\).  
\end{definition}

\begin{lemma}[\cite{JoyalTierney:1984}, \cite{BB:1992}]\label{Lem:5}
	For any frame \(L\) and any element \(a \in L\) there is a unique extension \(\map{x_a}{L}{L_a}\) in which \(a\) is freely complemented. 
\end{lemma}

\begin{proof}	
	Let \(x_a\) be the product of the closed and open quotient maps corresponding to \(a\).
	\[
	\begin{tikzcd}
		L			\arrow{r}{o_a}
					\arrow{d}[swap]{c_a}
					\arrow{dr}{x_a}
		&\downset{a}\\
		\upset{a}	
		&L_a \equiv \upset{a}\times\downset{a}
					\arrow{u}
					\arrow{l}
	\end{tikzcd},
	\]
	Consider an arbitrary extension \(\map{k}{L}{K}\) in which \(a\) is complemented, and define 
	\[
		\map{n}{L_a}{K}
		= \Big((b, c) \mapsto \big(k(b) \wedge k(a)^*\big) \vee k(c)\Big),
		\qquad (b, c) \in L_a.
	\] 
	\[
		\begin{tikzcd}
			L 			\arrow{r}{x_a}
						\arrow{d}[swap]{k}
			&L_a 
						\arrow{dl}{n}\\
			K
			&
		\end{tikzcd}
	\]
	Then routine calculation establishes that \(n\) is the unique frame homomorphism making the diagram commute.
\end{proof}

\begin{lemma}\label{Lem:7}
	Any frame \(L\) has a unique extension \(\map{p'}{L}{L'}\) in which every cozero element is freely complemented. 
\end{lemma}

\begin{proof}
	Let \(\{L_a \xra{y_a} L'\}_{\coz L}\) be the pushout of the source \(\{L \xra{x_a}L_a\}_{\coz L}\). The desired extension is \(p' \equiv y_a \circ x_a\) for some (any) \(a \in \coz L\).
\end{proof}

\subsection*{The \(P\)-frame reflection}

\begin{definition*}
	A frame \(L\) is called a \emph{\(P\)-frame} if each cozero element of \(L\) is complemented. We denote the full subcategory of \(\mbf{Frm}\) comprised of the \(P\)-frames by \(\mbf{PFrm}\).
\end{definition*}

\begin{theorem}[{\cite[7.13]{BallWWZenk:2010}}]
	\(\mbf{PFrm}\) is a monoreflective subcategory of \(\mbf{Frm}\).
\end{theorem}

\begin{proof}[Proof outline]
	Iterate the extension \(\map{p'}{L}{L'}\) of Lemma \ref{Lem:7} transfinitely. That the iteration stablizes is nontrivial, and hinges on the fact that the Lindel\"{o}f degree does not grow.
\end{proof}

Theorem \ref{Thm:4} clarifies the significance of \(P\)-frames for the theory of vector lattices.

\begin{definition*}[epicomplete \(\mbf{V}\)-object, \(\mbf{ecV}\)]
	A vector lattice \(G\) is called \emph{epicomplete} if it has no proper epic extensions in \(\mbf{V}\). We denote the full subcategory of \(\mbf{V}\) comprised of the epicomplete objects by \(\mbf{ecV}\).
\end{definition*}

\begin{theorem}[\cite{BallHager:1990}]\label{Thm:5}
	\(\mbf{ecV}\) is a monoreflective subcategory of \(\mbf{V}\).
\end{theorem}

\begin{theorem}[{\cite[4.4.2]{BallHagerWW:2015}}]\label{Thm:4}
	The epicomplete objects in \(\mbf{V}\) are those of the form \(\mcal{R}P\) for a \(P\)-frame \(P\). 
\end{theorem}

\begin{definition*}[\(C^*\)-quotient]
	A frame surjection \(\map{q}{L}{M}\) is called a \emph{\(C^*\)-quotient} if every bounded frame homomorphism \(\map{g}{\mcal{O}\mbb{R}}{M}\) factors through \(q\).
\end{definition*}

\section{A Yosida Theorem}\label{Sec:YosThm}

\subsection*{Compactifications and \(\mcal{E}q\)}

\begin{definition*}[compactification, \(\map{(l, m)}{q}{r}\), \(\mbf{Cmp}\), \(\map{p}{\mcal{O}\exR}{\mcal{O}\mbb{R}}\), domain of reality]	
	A \emph{compactification} is a dense surjective frame homomorphism $\map{q}{K}{M}$ with compact domain $K$. For compactifications $\map{q}{K}{M}$ and $\map{r}{L}{N}$, a \emph{homomorphism of compactifications} is a pair $(l, m)$ of frame homomorphisms such that $m \circ q = r \circ l$.  
	\[
	\begin{tikzcd}
		K 					\arrow{r}{l} 		
							\arrow{d}[swap]{q}
		& L 				\arrow{d}{r}\\ 
		M 					\arrow{r}[swap]{m}
		& N 
	\end{tikzcd}	
	\]	
	We write $\map{(l, m)}{q}{r}$. We denote the category of compactifications with their homomorphisms by $\mbf{Cmp}$. Of central importance is the compactification which is the frame homomorphism of the inclusion $\mbb{R} \to \exR$ of the real numbers \(\mbb{R}\) in the extended real numbers \(\ol{\mbb{R}} = \mbb{R} \cup \{\pm \infty\}\), given by the formula
	\[
		\map{p}{\mcal{O}\exR}{\mcal{O}\mbb{R}}
		=(U \mapsto U \cap (-, -)),
		\qquad U \in \mcal{O}\exR.
	\] 
	The \emph{domain of reality} of a frame homomorphism \(\map{h'}{\mcal{O}\exR}{L}\) is the cozero element \(h'(-,-) \in L\). (Here \((-,-)\) designates \(\setof{x \in \ol{\mbb{R}}}{-\infty < x < \infty} \in \mcal{O}\ol{\mbb{R}}\).)
\end{definition*}

\begin{lemma}\label{Lem:21}
	Let $\map{q}{K}{M}$ be a compactification. 
	\[
		\begin{tikzcd}
		\mcal{O}\exR 				\arrow{r}{h'} 
									\arrow{d}[swap]{p} 
		&K 							\arrow{d}{q}\\
		\mcal{O}\mbb{R} 			\arrow{r}[swap]{h} 
		&M
		\end{tikzcd}
	\]
	\begin{enumerate}
		\item (\cite[3.3.1]{Ball:2021}) 
		A homomorphism $h'$ drops to a homomorphism $h$ such that $h \circ p = q \circ h'$ if and only if \(q\) maps the domain of reality of \(h'\) to the top of \(M\), i.e., $q \circ h'(-, -) = \top$. 
		
		\item (\cite[5.4.1]{Ball:2021}) 
		A homomorphism \(h\) extends to a homomorphism \(h'\) such that $h \circ p = q \circ h'$ if and only if \(h \wedge n\) factors through \(q\) for all \(n \in \mbb{N}\).  
	\end{enumerate}
\end{lemma}

\begin{definition*}[\(\mcal{E}q\)]
	For a given compactification \(\map{q}{K}{M}\), we use \(\mcal{E}q\) to denote the family of compactification morphisms \(\map{(h', h)}{p}{q}\) which make the diagram of Lemma \ref{Lem:21} commute. In symbols,
	\[
		\mcal{E}q
		\equiv \setof{\map{(h', h)}{p}{q}}{q \circ h' = h \circ p}.
	\]
	\(\mcal{E}q\) inherits the (partial) vector lattice operations from \(\mcal{R}M\) as follows. Let \(\oplus\) designate any one of the operations \(+\), \(-\), \(\wedge\), or \(\vee\), and let \((h_i', h_i) \in \mcal{E}q\). Then we define \((h_1', h_1) \oplus (h_2', h_2) = (h_3', h_3)\) in \(\mcal{E}q\) to mean that \(h_1 \oplus h_2 = h_3\) in \(\mcal{R}M\).
\end{definition*}
		
		The issue that arises when attempting to add two elements \((h_i', h_i) \in \mcal{E}q\) is that the sum \(h_1 + h_2\), which always exists in \(\mcal{R}M\), may have no extension to a frame homomorphism \(\mcal{O}\ol{\mbb{R}} \to K\) which makes the diagram of Lemma \ref{Lem:21} commute. Lemma \ref{Lem:1} gives the details.   

\begin{definition*}[\(C^*\)-quotient]
	A frame surjection \(\map{q}{L}{M}\) is called a \emph{\(C^*\)-quotient} if every bounded frame homomorphism \(\map{g}{\mcal{O}\mbb{R}}{M}\) factors through \(q\).
\end{definition*}

\begin{lemma}[{\cite[6.1.1]{Ball:2021}}]\label{Lem:1}
	For any compactification \(\map{q}{K}{M}\), \(\mcal{E}q\) is closed under the operations of negation, meet, join, and scalar multiplication.  It is closed under addition, and is therefore a vector lattice, if and only if the open quotient of each cozero element \(a \in K\) such that \(q(a) = \top\) is a \(C^*\)-quotient. 
\end{lemma}

Example \ref{Ex:1} is an instance in which \(\mcal{E}q\) fails to be a vector lattice.   

\begin{example}\label{Ex:1}
	Let \(\map{q}{\mcal{O}[0, 1]}{\mcal{O}(0, 1]}\) be the frame map of the inclusion \((0, 1] \to [0, 1]\), and let \(f,g \in \mcal{R}(0, 1]\) be the frame maps of the continuous functions \(\tilde{f}, \tilde{g} \in \mcal{C}(0,1]\) defined by the rules \(\tilde{f}(x) = 1/x\) and \(\tilde{g} = 1/x + sin(1/x)\), \(x \in (0,1 ]\). Then both \(f\) and \(g\) satisfy the condition of Lemma \ref{Lem:21}(2) and therefore define members of \(\mcal{E}q\). But the same cannot be said of \(g - f\).
\end{example}

\begin{definition*}[vector lattice snugly embedded in \(\mcal{E}q\)]
	Let \(\map{q}{K}{M}\) be a compactification. If a subset \(G \subseteq \mcal{E}q\) is closed under all of the vector lattice operations we shall say that \emph{\(G\) is a vector lattice in \(\mcal{E}q\).} We shall say that a vector lattice \(G\) in \(\mcal{E} q\) is \emph{snugly embedded} if \(\coz G\) join-generates \(M\), and if \(q\) is the cointersection (pushout) of the open quotients of the domains of reality of the positive elements of \(G\). (Note that this condition implies that \(M\) is Lindel\"{o}f.)  
\end{definition*}

Proposition \ref{Prop:5} is an aside. It's proof is couched in terms of the discussion in Section \ref{Sec:RemLit}.

\begin{definition*}[majorizing subset, cofinality]
	A subset \(G_0\) of a vector lattice \(G\) is said to be \emph{majorizing} if 
	\[
		\forall g\in G^+\, \exists g_0 \in G_0\ (g \leq g_0).
	\]
	The \emph{cofinality} of a vector lattice \(G\) is the least regular cardinal number \(\kappa\) for which there is a majorizing subset \(G_0 \subseteq G^+\) of cardinality strictly less than \(\kappa\). In case \(\kappa = \omega_1\), we say that \(G\) has \emph{countable cofinality}.
\end{definition*}

\begin{proposition}\label{Prop:5}
	Let \(G\) be a vector lattice snugly embedded in \(\mcal{E}q\) for the compactification \(\map{q}{K}{M}\). If \(G\) has countable cofinality then \(M\) is spatial.
\end{proposition}

\begin{proof}
	Since \(K\) is compact it is spatial, i.e., the topology of a unique compact Hausdorff space \(X\) and, according to the classical theorem of Hager and Robertson (\cite{HagerRobertson:1977}), \(G\) is isomorphic to a vector lattice in \(\mcal{D}X\).  In this representation, the domains of reality of the positive elements of \(G\) are dense open subsets of \(X\). The result follows from the Baire Category Theorem.
\end{proof}

\subsection*{The adjunction}
Theorem \ref{Thm:3} is an amalgam of the generalized Yosida Theorem of Hager and Robertson \cite{HagerRobertson:1977} with the pointfree Yosida Theorem of Madden and Vermeer \cite{MaddenVermeer:1986}. See Section \ref{Sec:RemLit} for further remarks concerning its connections to the existing literature.

\begin{theorem}[{\cite[3.3.2]{Ball:2021}}]\label{Thm:3} 
	For every vector lattice \(G\) there is a compactification $q$, a snugly embedded vector lattice $\widehat{G}$ in \(\mcal{E}q\), and a vector lattice isomorphism $\map{\mu_G}{G}{\widehat{G}}$ with the following universal property. For any compactification \(r\) and any vector lattice homomorphism from \(G\) onto a vector lattice in \(\mcal{E}r\) there is a unique homomorphism \(\map{(l,m)}{q}{r}\) which makes the diagram commute.  
	\[
		\begin{tikzcd}
			G 			\arrow{r}{\mu_G} 
						\arrow{ddr}[swap]{\theta}
			& \mcal{E}q \arrow{dd}{\mcal{E}(l,m)}
			& K 		\arrow{rr}{q} 
						\arrow{dd}[swap]{l}
			&&M 		
			&\\ 
			&&&\mcal{O}\ol{\mbb{R}}	\arrow{rr}[pos=.3]{p}
									\arrow{ul}[swap]{\hat{g}'}
									\arrow{dl}{\theta(g)'}
			&&\mcal{O}\mbb{R} 	\arrow{ul}[swap]{\hat{g}} 
									\arrow {dl} {\theta (g)}\\
			&\mcal{E}r
			&L	 		\arrow{rr}[swap]{r}
			&& N		\arrow[leftarrow,crossing over]{uu}[pos=.3]{m}
			& 
		\end{tikzcd}	
	\]
\end{theorem} 

\begin{notation*}
	We refer to 
	\begin{itemize}
		\item 
		\(\map{q}{K}{M}\) as the \emph{standard compactification of \(G\)}, 
		
		\item 
		\(K\) as the \emph{Yosida frame of \(G\)},
		
		\item 
		\(M\) as the \emph{Madden frame of \(G\)},
		
		\item 
		\(\map{\mu_G}{G}{\widehat{G}}\subseteq \mcal{E}q\) as the (pointfree) Yosida representation of \(G\). We abuse this terminology to the extent of using \(\mu_G\) to also represent the induced map \(G \to \mcal{R}M = (g \mapsto \hat{g})\), \(g \in G\).
	\end{itemize}
\end{notation*}

\section{Pointwise convergence}\label{Sec:PtwiseCon}
Pointwise convergence has a natural and elegant formulation in the pointfree context, and it takes on greater importance when vector lattices are analyzed in that context. The topic was introduced in \cite{BallHagerWW:2015} and exploited in \cite{Ball:2018} and \cite{Ball:2021}. In connection with the definition of pointwise suprema below, the reader should recall that for elements \(f,g \in \mcal{R}M\), \(f \leq g\) if and only if \(f(q, -) \leq g(q, -)\) for all \(q \in \mbb{Q}\) if and only if \(f(-, q) \geq g(-, q)\) for all \(q \in \mbb{Q}\) (\cite[3.1.4]{BallHager:1991}).

\begin{definition*}[\(f_0 = \bigvee^\bullet F\), \(f_0 = \bigwedge^\bullet F \)]
	We say that an element \(f_0 \in \mcal{R}M\) is the \emph{pointwise join} (\emph{pointwise meet}) of a subset \(F \subseteq \mcal{R}M\) if \(\bigvee_F f(r, -) = f_0(r, -)\) \big(\(\bigvee_F f(-, r) = f_0(-, r)\)\big) for all \(r \in \mbb{R}\). (Actually, it is sufficient for this condition to hold only for all \(q \in \mbb{Q}\).) We write \(\bigvee^\bullet F = f_0\) \big(\(\bigwedge^\bullet F = f_0\)\big). We say that an element \(g_0 \in G \) is the pointwise join (pointwise meet) of a subset \(G_0 \subseteq G\) if \(\bigvee^\bullet \widehat{G}_0 = \hat{g}_0\) \big(\(\bigwedge^\bullet \widehat{G}_0 = \hat{g}_0\)\big).
\end{definition*}

Propositions 4.1.3 and 4.22--25 of \cite{BallHagerWW:2015} articulate some of the nice features of pointwise suprema and infima; we quote the last two as Proposition \ref{Prop:3} below. The second feature, context freeness, may be taken as a definition of pointwise suprema and infima, independent of the representation of \(G\) as a vector lattice in \(\mcal{R}M\).

\begin{proposition}\label{Prop:3}
	Let \(G_0\) be a subset and \(g_0\) an element of \(G\). 
	\begin{enumerate}
		\item 
		Pointwise joins and meets have \emph{countable character}. That is, if \(\bigvee_{G_0}^\bullet g = g_0\) then \(\bigvee_{G_1}^\bullet g = g_0\) for some countable subset \(G_1 \subseteq G_0\), and dually for pointwise meets.
		
		\item 
		Pointwise joins and meets are \emph{context free.} That is, \(\bigvee_{G_0}^\bullet g = g_0\) if and only if \(\bigvee_{G_0} \theta(g) = \theta(g_0)\) for any \(\mbf{V}\)-homomorphism out of \(G\), and dually for pointwise meets. 
	\end{enumerate}
\end{proposition}

The notion of pointwise suprema and infima leads immediately to the appropriate notion of directional pointwise convergence. Sections 5 and 6 of \cite{Ball:2018} constitute a thorough development of this idea.

\begin{definition*}[\(g_n \searrow g_0\), \(g_n \nearrow g_0\)]
	A sequence \(\{g_n\} \subseteq G\) \emph{converges downwards (upwards) to a limit \(g_0\)} if it is decreasing (increasing), i.e., \(g_{n + 1} \leq g_n\) (\(g_{n + 1} \geq g_n\)) for all \(n\), and \(\bigwedge_n^\bullet g_n = g_0\) (\(\bigvee_n^\bullet g_n = g_0\)). We write \(g_n \searrow g_o\) (\(g_n \nearrow g_0\)). 
\end{definition*}

Unqualified pointwise convergence is obtained from directional pointwise convergence by means of the lim sup definition, that is, by requiring the pointwise downwards convergence of the suprema of the tails of the
sequence. This makes sense if the context assures the existence of the aforementioned suprema. It turns out that the appropriate context is the epicompletion of \(G\), which for our purposes is realized by the extension
\(G \xra{\mu_G} \mcal{R}M \xra{\mcal{R}p} \mcal{R}P\), where \(\map{p}{M}{P}\) is the \(P\)-frame reflection of \(M\) introduced in Theorems \ref{Thm:5} and \ref{Thm:4}. Lemmas \ref{Lem:9} implies  Lemma \ref{Lem:8}; together they provide evidence of the utility of the extension \(G \to \mcal{R}P\) for our purposes.

\begin{lemma}[{\cite[3.9]{BallHager:1990a}}]\label{Lem:9}
	A vector lattice is of the form \(\mcal{R}P\) for some \(P\)-frame \(P\) if and only if it is countably conditionally and laterally complete. 
\end{lemma}

\begin{lemma}[{\cite[4.6]{Ball:2018}}]\label{Lem:8}
	When \(\mcal{R}M\) is viewed as a subobject of \(\mcal{R}P\), any sequence in \(\mcal{R}M\) which is bounded above (below) in \(\mcal{R}P\) has a pointwise join (meet) in \(\mcal{R}P\). 
\end{lemma}

\begin{definition*}[pointwise convergence]
	We say that a sequence $\{g_n\} \subseteq G^+$ \emph{converges pointwise to $0$}, and write $g_n \pc 0$, if 
	\[
		\sbv{r}{j \geq n}^\bullet g_j \searrow 0 
	\]
	holds in $\mcal{R}_0 P$, i.e., if $\bigwedge_n^\bullet \bigvee_{j \geq n}^\bullet g_j = 0$ in $\mcal{R} P$. An arbitrary sequence $\{g_n\} \subseteq G$ converges pointwise to an element $g_0$ if $|g_n - g_0| \pc 0$. 
\end{definition*}

We reiterate for emphasis that the joins that appear in the definition of pointwise convergence need not exist in $\mcal{R} M$, but do exist in $\mcal{R}P$. 

\begin{proposition}\label{Prop:6}
	Let $\{g_n\}$ be a bounded sequence in $G^+$. Then $g_n \pc 0$ if and only if 
	\[
		\forall \epsilon > 0\ \Bigg(\bigvee_m \sbw{r}{n \geq m} g_n(-\infty,  \epsilon) = \top \Bigg)
	\]
	holds in $P$. 
\end{proposition}

\begin{proof}
	Let $h_m \equiv \bigvee_{n \geq m}^\bullet g_n$. To say that $h_m \searrow 0$ is to say that $\bigwedge_m^\bullet h_m = 0$, i.e., that $\bigvee_m h_m(-\infty, \epsilon) = \top$ for all $\epsilon > 0$. Now it is a general fact that $f(-\infty, \epsilon) = \bigvee_{\delta < \epsilon}f(\delta, \infty)^*$ for any $f \in \mcal{R}_0 L$ and $\epsilon \in \mbb{R}$ (\cite[3.1.1]{BallHager:1991}), hence
	\begin{align*}
		h_m(-\infty, \epsilon) 
		&= \Bigg(\sbv{r}{n \geq m}^\bullet g_n\Bigg)(-\infty, \epsilon)
		= \sbv{l}{\delta < \epsilon}\Bigg(\sbv{r}{n \geq m}^\bullet g_n\Bigg)(\delta, \infty)^*\\
		&= \sbv{l}{\delta < \epsilon}\Bigg(\sbv{r}{n \geq m} g_n(\delta, \infty)\Bigg)^*
		= \sbv{l}{\delta < \epsilon}\sbw{r}{n \geq m} g_n(\delta, \infty)^*.
	\end{align*}
	We have 
	\begin{align*}
		g_n \pc 0 
		&\iff h_m \searrow 0
		\iff \forall \epsilon  > 0\ \Bigg(\bigvee_m \bigvee_{\delta <  \epsilon}\sbw{r}{n \geq m}g_n(\delta, \infty)^* = \top\Bigg)\\
		&\iff \forall \epsilon  > 0\ \Bigg(\bigvee_{\delta <  \epsilon}\bigvee_m \sbw{r}{n \geq m}g_n(\delta, \infty)^* = \top\Bigg)
	\end{align*}	
	Since $\bigvee_m\bigwedge_{n \geq m}g_n(\delta, \infty)^*$ is an increasing function of $\delta$, 
	\[
		g_n \pc 0  
		\iff \forall \epsilon > 0\ \Bigg(\bigvee_m \sbw{r}{n \geq  m}g_n(\epsilon, \infty)^* = \top\Bigg). \qedhere
	\]
\end{proof}

Proposition \ref{Prop:16} points out that pointwise convergence generalizes directional pointwise convergence. 

\begin{proposition}\label{Prop:16}
	For a bounded sequence $\{g_n\} \subseteq G^+$, 
	\[
		g_n \searrow 0 
		\iff \text{$\{g_n\}$ is decreasing and $g_n \pc 0$},
	\]
	and dually.
\end{proposition}

\begin{proof}
	If $g_n \searrow 0$ then $\{g_n\}$ is decreasing and $\bigwedge_n^\bullet g_n = 0$, hence $\bigvee_{j \geq n}g_j = g_n$ for all $n$, and $\bigwedge_n^\bullet \bigvee_{j \geq n}^\bullet g_j = \bigwedge_n^\bullet g_n = 0$, i.e., $g_n \pc 0$. On the other hand, suppose $\{g_n\}$ is decreasing and $g_n \pc 0$, i.e., $\bigwedge_n^\bullet \bigvee_{j \geq n}g_j^\bullet = 0$. Since once again $\bigvee_n^\bullet = g_n$, this says directly that $\bigwedge_n^\bullet g_n = 0$, i.e., that $g_n \searrow 0$.
\end{proof}

\section{Relative uniform convergence}

\begin{definition*}[regulator, \(g_n \rightarrow g_0\, (f)\), relative uniform convergence \(g_n \xrightarrow{u} g_0\)]
	A sequence \(\{g_n\} \subseteq G\) is said to \emph{converge to an element \(g_0\) with regulator \(f \in G^+\)}, written \(g_n \rightarrow g_0\, (f)\), if for every positive integer \(k\) there exists an integer \(m_k\) such that \(\abs{g_0 - g_n} \leq f/k\) for all \(n \geq m_k\). In symbols,  
	\[
		\forall k \, \exists m_k\, \forall n \geq m_k\ \big(\abs{g_0 - g_n} \leq f/k \big).
	\]
	The sequence is said to \emph{converge relatively uniformly to \(g_0\)}, written \(g_n \xrightarrow{ru} g_0\) if there exists an element \(f \in G^+\) for which \(g_n \rightarrow g_0\, (f)\).  
\end{definition*}

We record our notation for open subsets of the reals, and then provide a few preliminary results which will facilitate the material that follows.  

\begin{notation*}[\(r^\varepsilon\), \(U \subseteq_\varepsilon V\), \(U^\circ_\varepsilon\)]
	For a real number \(r\) and positive real number \(\varepsilon\), we denote the \(\varepsilon\)-neighborhood of \(r\) by 
	\[
		r^\varepsilon
		\equiv \setof{s \in \mbb{R}}{ \abs{r - s} < \varepsilon}
	\]
	When speaking of subsets \(U,V \subseteq \mbb{R}\), we say that  \emph{\(U\) is \(\varepsilon\)-contained in \(V\)}, and write \(U \subseteq_\varepsilon V\), if \(u^\varepsilon \subseteq V\) for all \(u \in U\). The \emph{\(\varepsilon\)-interior of \(V\)} is the largest open subset of \(\mbb{R}\) which is \(\varepsilon\)-contained in \(V\), i.e.,
	\[
		V^\circ_\varepsilon
		\equiv \bigcup\setof{U \in \mcal{O}\mbb{R}}{U \subseteq_\varepsilon  V}
		= \setof{v}{v^\varepsilon \subseteq V}.
	\]
\end{notation*}

\begin{lemma}\label{Lem:17}
	For all \(f,g \in \mcal{R} M\) and all \(\varepsilon > 0\),
	\begin{align*}
		\abs{f - g}(\varepsilon, -)
		&= \bigvee_s \big((f(-, s) \wedge g(s + \varepsilon, -)) \vee (g(-, s) \wedge f(s + \varepsilon, -))\big),\ \text{ and}\\
		\abs{f - g}(-, \varepsilon)
		&= \bigvee_s \left(f(s^{\varepsilon/2}) \wedge g(s^{\varepsilon/2})\right).
	\end{align*}
\end{lemma}

\begin{proof}
	According to \cite[3.1.1]{BallWalters:2002},
	\[
	\abs{f - g}(\varepsilon, -)
	= \sbv{lr}{\abs{U_f - U_g}\subseteq (\varepsilon, -)}\left(f(U_f) \wedge g(U_g)\right),
	\]
	where \(U_f\) and \(U_g\) range over open subsets of \(\mbb{R}\) such that \(\abs{u - v} > \varepsilon\) for all \(u \in U_f\) and \(v \in U_g\).  Thus either \(U_f\) is bounded above and \(U_g\) is bounded below, or vice-versa.  In the first case there exists a real number \(s\) such that \(U_f \subseteq (-, s)\) and \(U_g \subseteq (s, -)\), and in the second case there exists a real number \(s\) such that \(U_f \subseteq (s, -)\) and \(U_g \subseteq (-, s)\).  
	
	To verify the second displayed identity, consider open sets \(U_f, U_g \in \mcal{O}\mbb{R}\) such that \(\abs{U_f - U_g} \subseteq (-, \varepsilon)\). If \(u \in U_f\) and \(v \in U_g\) then the midpoint \(s \equiv (u + v)/2\) has the feature that \(u, v \in s^{\varepsilon/2}\) and \(\abs{s^{\varepsilon/2} - s^{\varepsilon/2}} \subseteq (-, \varepsilon)\). 
\end{proof}

\begin{lemma}\label{Lem:19}
	For any \(f,g \in \mcal{R} M\), any \(U \in \mcal{O}\mbb{R}\), and any \(\varepsilon > 0\),
	\[
		\abs{f - g}(-, \varepsilon) \wedge f(U^\circ_\varepsilon)
		\leq g(U).
	\]
\end{lemma}

\begin{proof}
	From Lemma \ref{Lem:17} we get
	\[
	\abs{f - g}(-, \varepsilon) \wedge f(U^\circ_\varepsilon)
	= \bigvee_s \left(f(U^\circ_\varepsilon) \wedge f(s^{\varepsilon/2}) \wedge g(s^{\varepsilon/2})\right)
	= \bigvee_s \left(f(U^\circ_\varepsilon \wedge s^{\varepsilon/2}) \wedge g(s^{\varepsilon/2})\right).
	\]
	We claim that \(s^{\varepsilon/2} \subseteq U\) for any real number \(s\) such that \(U^\circ_\varepsilon \wedge s^{\varepsilon/2} \neq \emptyset\). For if \(r \in s^{\varepsilon/2}\) and \(t \in U^\circ_\varepsilon \wedge s^{\varepsilon/2}\) then \(\abs{r - t} = \abs{r - s} + \abs{s - t} < \varepsilon/2 + \varepsilon/2 = \varepsilon\), so that \(r \in U\) by virtue of the fact that \(U^\circ_\varepsilon \subseteq_\varepsilon U\). In light of the equation displayed above, the claim proves the lemma.
\end{proof}

We supply a couple of alternative characterizations of indicated relative uniform convergence, and then list some of its nice properties.

\begin{lemma}\label{Lem:16}
	For a sequence \(\{g_n\} \subseteq G^+\) and element \(0 < f, g_0 \in G\), the first three conditions are equivalent to one another.  From this  follows the equivalence of the second three conditions to one another.  
	\begin{enumerate}
		\item 
		\(g_n \to 0\, (f)\).
		
		\item 
		For every integer \(k\) there exists an integer \(m_k\) such that for all \(n \geq m_k\) and all \(\varepsilon > 0\),
		\[
			g_n(0^\varepsilon) \geq f(0^{k\varepsilon}).
		\] 
		
		\item 
		For every integer \(k\) there exists an integer \(m_k\) such that for all \(n \geq m_k\) and all \(\varepsilon > 0\),
		\[
			g_n(\varepsilon, -) \leq f(k\varepsilon, -).
		\]

		\item 
		\(g_n \rightarrow g_0\, (f)\).	
		
		\item 
		For every integer \(k\) there exists an integer \(m_k\) such that for all \(n \geq m_k\) and all \(\varepsilon > 0\),
		\[
			\bigvee_t \big(g_n(t^{\varepsilon/2}) \wedge  g_0(t^{\varepsilon/2})\big) \geq f(0^{k\varepsilon}).	
		\]
		
		\item 
		For every integer \(k\) there exists an integer \(m_k\) such that for all \(n \geq m_k\) and all \(\varepsilon > 0\),
		\[
			\bigvee_t \Big(\big(g_0(t + \varepsilon, -) \wedge g_n(-, t)\big) \vee \big(g_n(t + \varepsilon, -) \wedge g_0(-, t)\big)\Big) \leq f(k\varepsilon, -).
		\]
	\end{enumerate}
\end{lemma}

\begin{proof}
	To prove the equivalence of the first three conditions simply observe that \(0 \leq g_n \leq f/k\) if and only if \(g_n(\varepsilon, -) \leq (f/k)(\varepsilon, -) = f(k\varepsilon, -)\) for all \(\varepsilon > 0\) if and only if \(g_n(-, \varepsilon) \geq (f/k)(-, \varepsilon) = f(-, k\varepsilon)\) for all \(\varepsilon > 0\). To prove the equivalence of conditions (4) and (6) use the first equation displayed in Lemma \ref{Lem:17} to express  \(\abs{g_n - g_0}(\varepsilon, -)\) as 
	\[
		\bigvee_t\Big(\big(g_0(t + \varepsilon, -) \wedge g_n(-, t)\big) \vee \big(g_n(t + \varepsilon, -) \wedge g_0(-, t)\big) \Big).
	\]
	To prove the equivalence of conditions (4) and (5) use the second equation of the lemma in similar fashion. 
\end{proof}

\begin{proposition}[{\cite[14]{BallHager:1999b}}]\label{Prop:22}
	Vector lattice homomorphisms are continuous with respect to relative uniform convergence. That is, for any \(\mbf{V}\)-morphism $\tau \colon G \to K$ and any sequence $\{g_n\}$ and element $g_0$ in $G$,
	\[
		g_n \xrightarrow{ru} g_0 
		\implies \tau(g_n) \xrightarrow{ru} \tau(g_0).
	\]
	Furthermore, the converse holds if \(\tau\) is one-one, i.e., if \(\tau(g_n) \to \tau(g_0)\,(\tau(f))\) for some sequence \(\{g_n\} \subseteq G\) and elements \(f, g_0 \in G\) with \(f > 0\) then \(g_n \to g_0\,(f)\). 
\end{proposition}

\begin{proof}
	Suppose \(g_n \to g_0\, (f)\), let $\mcal{K} K \equiv N$ and let $m \colon M \to N$ be the pointed frame homomorphism which realizes $\tau$, i.e., $\tau(g) = k \circ g$ for all $g \in G$. Our asumption is that
	\[
		\forall k \, \exists m \, \forall n \geq m\ (\abs{g_0 - g_n} \leq f/k).
	\]
	An application of the \(\mbf{V}\) homomorphism \(k\) yields \(\tau(g_n) \rightarrow \tau(g_0)\, (\tau(f))\), i.e., \(\tau(g_n) \xrightarrow{ru} \tau(g_0)\).
\end{proof}

The vector lattice operations are relatively uniformly continuous. 

\begin{proposition}\label{Prop:23}
	Let $\{f_n\}$ and $\{g_n\}$ be sequences in $G$, and let $f_0$ and $g_0$ be elements of $G$.
	\begin{enumerate}
		\item 
		$\dot{g}_0 \xrightarrow{ru} g_0$, where $\dot{g}_0$ designates the constant sequence, i.e., $(\dot{g}_0)_n = g_0$ for all $n$.
		
		\item 
		$g_n \xrightarrow{ru} g_0$ if and only if $(-g_n) \xrightarrow{ru} (-g_0)$.
		
		\item
		If $f_n \geq g_n \geq 0$ and $f_n \xrightarrow{ru} 0$ then $g_n \xrightarrow{ru} 0$.
			
		\item 
		If $g_n \xrightarrow{ru} g_0$ and $r \in \mbb{R}$ then $(rg_n) \xrightarrow{ru} rg_0$.

		\item 
		If $f_n \xrightarrow{ru} f_0$ and $g_n \xrightarrow{ru} g_0$ then $(f_n	\oplus g_n) \xrightarrow{ru} (f_0 \oplus g_0)$, where $\oplus$ stands for one of the operations $+$, $-$, $\vee$, or $\wedge$. 
		
		\item 
		If $g_n \xrightarrow{ru} g_0$ and $g_n \xrightarrow{ru} f_0$ then $g_0 = f_0$.
	\end{enumerate}
\end{proposition}

\begin{proof}
	The verifications are straightforward. 
\end{proof}

\begin{proposition}\label{Prop:24}
	Pointwise convergence is finer than relative uniform convergence.  That is, for any sequence \(\{g_n\} \subseteq G\) and any element \(g_0 \in G\), if \(g_n \xrightarrow{ru} g_0\) then \(g_n \pc g_0\).
\end{proposition}

\begin{proof}
	We treat the case in which \(g_n \rightarrow g_0\, (f)\), for which it is sufficient to handle the subcase in which \(\{g_n\} \subseteq G^+\), \(g_0 = 0\), and \(f \in G^+\). According to Proposition \ref{Prop:6}, we must show that \(\bigvee_m \bigwedge_{n \geq m}g_n(0^\varepsilon) = \top\) for all \(\varepsilon > 0\). For that purpose fix \(\varepsilon > 0\); for each integer \(l\) let \(k\) be the least integer greater than \(l/\varepsilon\), and use Lemma \ref{Lem:16}(2) to find an integer \(m_k\) such that \(g_n(0^\varepsilon) \geq f(0^{k\varepsilon})\) for all \(n \geq m_l\). Then \(\bigwedge_{n \geq m_k}g_n(0^\varepsilon) \geq f(0^{k\varepsilon}) \geq f(0^l)\). 
\end{proof}

Corollary \ref{Cor:1} will be put to use in the proof of Proposition \ref{Prop:4}.

\begin{corollary}\label{Cor:1}
	If \(g_n \xra{ru} g_0\) for some increasing sequence \(\{g_n\} \subseteq G^+\) then \(g_n \nearrow g_0\), with the consequence that \(\coz g_0 = \bigvee_n \coz g_n\).
\end{corollary}

\begin{proof}
	It is Proposition \ref{Prop:16} which informs us \(g_n \nearrow g_0\), which is to say that \(\bigvee_n^\bullet g_n = g_0\). The meaning of the latter condition is that \(\bigvee_n g_n(r, -) = g_0(r, -)\) for all \(r \in \mbb{R}\). Letting \(r = 0\) gives the desired result. 
\end{proof}

\subsection*{Inner and outer relative uniform closure operators}

Relative uniform convergence gives rise to two distinct closure operators.

\begin{definition*}[the outer (inner) relative uniform completion of \(G\) in \(H\), \(\oruc_H G\), (\(\iruc_H G\)), oru-complete (iru-complete)]
	Suppose that \(G\) is a sub-vector lattice of \(H\). Let
	\begin{align*}		
		G^1_H
		&\equiv \setof{h \in H}{\exists \, \{g_n\} \subseteq G\, \exists f  \in H^+\ \big(g_n \to h\, (f)\big)}\\
		\Big(\prescript{1}{}{G}_H
		&\equiv \setof{h \in H}{\exists \, \{g_n\} \subseteq G\, \exists f  \in G^+\ \big(g_n \to h\, (f)\big)} \Big).
	\end{align*}
	We define the \emph{outer (inner) relative uniform closure of \(G\) in \(H\)} inductively as follows. For all ordinals \(\alpha\) define
	\begin{align*}
		G^0_H
		&\equiv G
		& \big(\prescript{0}{}{G}_H
		&\equiv G\big),
		& &\\
		G^\alpha_H
		&\equiv (G^\beta_H)^1_H
		& \Big(\prescript{\alpha}{}{G}_H
		&\equiv \prescript{1}{}{(\prescript{\beta}{}{G}_H})_H\Big)
		& &\text{for \(\alpha = \beta + 1\),}\\
		G^\alpha_H
		&\equiv \bigcup_{\beta < \alpha}G^\beta_H
		&\bigg(\prescript{\alpha}{}{G}_H
		&\equiv \bigcup_{\beta < \alpha}\prescript{\beta}{}{G}_H\bigg) 
		&&\text{for \(\alpha\) a limit ordinal,}\\
		\oruc_H G
		&\equiv G^\alpha_H
		& \Big(\iruc_H G
		&\equiv \prescript{\alpha}{}{G}_H\Big)
		&&\text{for \(\alpha\) such that \(G^\alpha_H = G^{\alpha + 1}_H\) \(\big(\prescript{\alpha}{}{G}_H = \prescript{\alpha + 1}{}{G}_H\big)\).}
	\end{align*}
	If \(\oruc_H G = H\) (\(\iruc_H G = H\)) we say that \emph{\(G\) is oru-dense (iru-dense) in \(H\)}. We say that \(G\) is \emph{oru-complete (iru-complete)} if \(\oruc_H G = G\) (\(\iruc_H G = G\)) for every vector lattice \(H\) which contains \(G\) as a sub-vector lattice. We denote the full subcategory of \(\mbf{V}\) comprised of the oru-complete (iru-complete) vector lattices by \(\mbf{orucV}\) (\(\mbf{irucV}\)).  
\end{definition*}

It is well known, and explicit in the development in \cite{BallHager:1999b}, that the iteration in the definition of \(\oruc_H G\) may proceed through \(\omega_1\) steps. It is less well known that the iteration in the definition of \(\iruc_H G\) may take more than one step; an example requiring two steps can be found in \cite{BallHager:2025}. We know of no example requiring more than two steps, but we think they exist. 

\begin{proposition}[{\cite[14]{BallHager:1999b}}]\label{Prop:1}
	If \(G\) is an oru-dense sub-vector lattice of \(H\), i.e., if \(\oruc_H G = H\), then \(H\) is an epic extension of \(G\).
\end{proposition}

\begin{proof}
	Consider vector lattice homomorphisms \(\map{\theta_i}{H}{K}\) which agree on \(G\). Then for any \(h \in H^+\) there is a sequence \(\{g_n\} \subseteq G\) and regulator \(f \in H^+\) such that \(g_n \to h\,(f)\). By Proposition \ref{Prop:22} it follows that \(\theta_i(g_n) \to \theta_i(h)\,(\theta_i(f))\) for \(i = 1,2\). But then the same sequence, labeled either \(\{\theta_1(g_n)\}\) or \(\{\theta_2(g_n)\}\), converges uniformly to both \(\theta_1(h)\) and \(\theta_2(h)\) with the same regulator \(\theta_1(f) \vee \theta_2(f))\). Since uniform limits are unique, we conclude that \(\theta_1(h) = \theta_2(h)\).
\end{proof}

\section{\(\mcal{R} M\) is ru-complete}\label{Sec:RMruc}

\begin{definition*}[ru-Cauchy sequence, ru-Cauchy complete vector lattice]
	A sequence \(\{g_n\}\) in a vector lattice \(G\) is said to be \emph{relatively uniformly Cauchy,} or simply \emph{ru-Cauchy,} if there is some regulator \(f \in G^+\) for which we have 
	\[
		\forall k\, \exists m_k\, \forall i,j \geq m_k\ \big(\abs{g_i - g_j} \leq f/k\big).
	\]
	\(G\) is said to be \emph{ru-Cauchy complete} if every ru-Cauchy sequence in \(G\) converges relatively uniformly to an element of \(G\).
\end{definition*}

Lemma \ref{Lem:13} provides a characterization of ru-Cauchy sequences in \(\mcal{R} M\).

\begin{lemma}\label{Lem:13}
	The following are equivalent for a sequence \(\{g_n\} \subseteq \mcal{R}^+ M\) and element \(f \in \mcal{R}^+M\).
	\begin{enumerate}
		\item 
		\(\{g_n\}\) is ru-Cauchy with regulator \(f\). 
		
		\item 
		For all \(k\) there exists an \(m_k\) such that for all \(i,j \geq m_k\) and all \(\varepsilon > 0\) we have
		\[
			\bigvee_s\Big(\big(g_i(-, s) \wedge g_j(s + \varepsilon, -)\big) \vee \big(g_j(-, s) \wedge g_i(s + \varepsilon, -)\big)\Big) 
			\leq f(k\varepsilon, -).
		\]
		
		\item 
		For all \(k\) there exists an \(m_k\) such that for all \(i,j \geq m_k\) and all \(\varepsilon > 0\) we have
		\[
			\bigvee_s \left(g_i(s^{\varepsilon/2}) \wedge g_j(s^{\varepsilon/2})\right)
			\geq f(-, k\varepsilon). 
		\]
	\end{enumerate} 
\end{lemma}

\begin{proof}
	By Lemma 3.1.4 of \cite{BallHager:1991} we get that \(\abs{g_i - g_j} \leq f/k\) if and only if 	
	\[
		\forall \varepsilon > 0\ \big(\abs{g_i - g_j}(\varepsilon, -) \leq (f/k)(\varepsilon, -)\big)
		\iff \forall \varepsilon > 0\ \big(\abs{g_i - g_j}(-, \varepsilon) \geq (f/k)(-, \varepsilon)\big). 
	\]
	Since \((f/k)(\varepsilon, -) = f(k\varepsilon, -)\), the first displayed identity of Lemma \ref{Lem:17} applied to the left condition above yields condition (2) of the lemma. Likewise the second displayed identity of Lemma \ref{Lem:17} applied to the right condition above yields condition (3) of the lemma. 
\end{proof}

\begin{lemma}\label{Lem:18}
	Let \(\{g_n\}\) be an ru-Cauchy sequence in \(\mcal{R} M\) with regulator \(f \in \mcal{R}^+_0 M\), and for each integer \(k\) let \(m_k\) be the least integer such that \(\abs{g_i - g_j} \leq f/k\) for all \(i,j \geq m_k\). Then for all \(k \in \mbb{N}\), \(r \in \mbb{R}\), and \(\varepsilon > 0\) we have
	\[
		g_{m_k}(-, r - \varepsilon) \wedge f(-, k\varepsilon) 
		\leq g_j (-, r)
		\qtq{and} g_{m_k}(r + \varepsilon, -) \wedge f(-, k\varepsilon) 
		\leq g_j (r, -)
	\] 
	for all \(j \geq m_k\).
\end{lemma}

\begin{proof}
	To check the left inequality, note that 
	\begin{align*}
		g_{m_k}(-, r - \varepsilon) \wedge f(-, k\varepsilon)
		&\leq g_{m_k}(-, r - \varepsilon) \wedge \abs{g_{m_k} - g_j}
		\leq g_j(-, r).
	\end{align*}
	Here the first inequality follows from Lemma \ref{Lem:13} and its proof, and the second is an application of Lemma \ref{Lem:19}. The argument for the right inequality is similar.
\end{proof}

\begin{theorem}[\cite{Hager:2015}]\label{Thm:12}
	\(\mcal{R} M\) is ru-Cauchy complete.
\end{theorem}

\begin{proof}
	Let \(\{g_n\}\) be an ru-Cauchy sequence in \(\mcal{R}^+ M\) with regulator \(f \in \mcal{R}^+_0 M\), and for each integer \(k\) let \(m_k\) be the least integer such that \(\abs{g_i - g_j} \leq f/k\) for all \(i,j \geq m_k\).  If the relative uniform limit of this sequence exists as an element \(g_0 \in \mcal{R} L\) then, since \(g_{m_k} - f/k \leq g_j \leq g_{m_k} + f/k\) for all \(j \geq m_k\), it must be the case that \(g_{m_k} - f/k \leq g_0 \leq g_{m_k} + f/k\) for all \(k\), so that \(g_{m_k} + f/k \to g_0\, (f)\) and hence \(g_{m_k} + f/k \pc g_0\) by Proposition \ref{Prop:24}. In other words, the only candidate for the limit is \(\bigwedge^\bullet_k(g_{m_k} + f/k)\); the only issue is whether this pointwise meet exists in \(\mcal{R} M\). 
	
	Thus we are led to define the function \(\map{g_0}{\setof{(-, r)}{r \in \mbb{R}}}{M}\) by the formula 
	\[
		g_0(-, r)
		\equiv \bigvee_k \big(g_{m_k} +f/k\big)(-, r),
		\qquad r \in \mbb{R}.
	\]
	Now 
	\begin{multline*}
		\big(g_{m_k} + f/k\big)(-, r)
		= \sbv{lr}{U + V \subseteq (-, r)}\big(g_{m_k}(U) \wedge (f/k)(V)\big)
		=  \sbv{lr}{\varepsilon > 0}\big(g_{m_k}(-, r - \varepsilon) \wedge (f/k)(-, \varepsilon)\big)\\
		= \sbv{lr}{\varepsilon > 0}\big(g_{m_k}(-, r - \varepsilon) \wedge f(-, k\varepsilon)\big),
	\end{multline*}
	so \(g_0(-, r) = \bigvee_k \bigvee_{\varepsilon > 0}\big(g_{m_k}(-, r - \varepsilon) \wedge f(-, k\varepsilon) \big)\). In order to show that \(g_0\) extends uniquely to a member of \(\mcal{R} L\), it is sufficient to establish that it has these three properties (\cite[3.1.2]{BallHager:1991}).
	\begin{enumerate}
		\item 
		\(\bigvee_{t > 0} g_0(-, r - t) = g_0(-, r)\) for all \(r \in \mbb{R}\).
		
		\item 
		\(g_0(-, r) \prec g_0(-, t)\) for all \(r < t\) in \(\mbb{R}\).
		
		\item 
		\(\bigvee_r g_0(-, r) = \bigvee_r g_0(-, r)^* = \top\).
	\end{enumerate}
	To verify (1) simply observe that 
	\begin{multline*}
		\sbv{lr}{t > 0}g_0(-, r - t)
		= \sbv{l}{t > 0}\bigvee_k\sbv{r}{\varepsilon > 0}\big(g_{m_k}(-, r - t - \varepsilon) \wedge f(-, k\varepsilon)\big)
		= \bigvee_k \sbv{l}{\varepsilon > 0} \sbv{r}{t > 0}\big(g_{m_k}-, r - t - \varepsilon) \wedge f(-, k\varepsilon)\big)\\
		= \bigvee_k \sbv{lr}{\varepsilon > 0} \left(f(-, k\varepsilon) \wedge \sbv{lr}{t > 0} g_{m_k}(r - t - \varepsilon)\right)
		= \bigvee_k \sbv{lr}{\varepsilon > 0} \left(f(-, k\varepsilon) \wedge g_{m_k}(-, r - \varepsilon)\right)
		= g_0(-, r).
	\end{multline*}
	
	To verify (2) fix real numbers \(r < t\) and put 
	\[
		a_r 
		\equiv \bigvee_l \sbv{lr}{\delta > 0}\big(g_{m_l}(r + \delta, -) \wedge f(-, l\delta)\big)
	\]
	We claim that \(a_r\) witnesses the fact that \(g_0(-, r) \prec g_0(-, t)\). First of all, 
	\begin{multline*}
		g_0(-, r) \wedge a_r
		= \bigvee_k \sbv{lr}{\varepsilon > 0}\big(g_{m_k}(- , r - \varepsilon) \wedge f(-, k\varepsilon) \wedge a_r\big)\\
		= \bigvee_k \sbv{lr}{\varepsilon > 0}\bigvee_l \sbv{lr}{\delta > 0}\left(g_{m_k}(- , r - \varepsilon) \wedge f(-, k\varepsilon) \wedge g_{m_l}(r + \delta, -) \wedge f(-, l\delta)\right).
	\end{multline*}
	For any integer \(k\) and \(\varepsilon > 0\) we have by Lemma \ref{Lem:18} that \(g_{m_k}(- , r - \varepsilon) \wedge f(-, k\varepsilon) \leq g_j(-, r)\) for all \(j \geq m_k\); likewise for any integer \(l\) and \(\delta > 0\) we have \( g_{m_l}(r + \delta, -) \wedge f(-, l\delta) \leq g_i(r, -)\) for all \(i \geq m_l\). Thus for any \(k,l \in \mbb{N}\) and \(\varepsilon, \delta > 0\) the term displayed in the quadruply indexed join above is bounded above by \(g_j(-, r) \wedge g_i(r, -)\) for any \(j > m_k\) and \(i > m_l\). By choosing \(i = j > m_k \vee m_l\) we see that this term, and therefore also the join, works out to \(\bot\).
		
	On the other hand we have
	\begin{multline*}
		g_0(-, t) \vee a_r
		= \bigvee_k \sbv{lr}{\varepsilon > 0}\big(g_{m_k}(- , t - \varepsilon) \wedge f(-, k\varepsilon) \vee a_r\big)\\
		= \bigvee_k \sbv{lr}{\varepsilon > 0}\bigvee_l \sbv{lr}{\delta > 0}\left(\big(g_{m_k}(- , t - \varepsilon) \wedge f(-, k\varepsilon)\big) \vee \big(g_{m_l}(r + \delta, -) \wedge f(-, l\delta)\big)\right).
	\end{multline*}
	When \(\varepsilon = \delta = (t - r)/3\) and \(k = l\), the term in the quadruply indexed join immediately above works out to \(f(-, k\varepsilon)\). Since \(\bigvee_k f(-, k\varepsilon) = \top\), we conclude that \(g_0(-, t) \vee a_r = \top\). The claim is proven.
	
	To verify (3) first note that 
	\begin{multline*}
		\bigvee_r g_0(-, r)
		= \bigvee_r \bigvee_k \sbv{lr}{\varepsilon > 0} \big(g_{m_k}(-, r - \varepsilon) \wedge f(-, k\varepsilon)\big)
		= \bigvee_k \sbv{lr}{\varepsilon > 0} \bigvee_r \big(g_{m_k}(-, r - \varepsilon) \wedge f(-, k\varepsilon)\big)\\
		= \bigvee_k \sbv{lr}{\varepsilon > 0} \left(f(-, k\varepsilon) \wedge \bigvee_r g_{m_k}(-, r - \varepsilon)\right)
		= \sbv{lr}{\varepsilon > 0} \bigvee_k f(-, k\varepsilon)
		= \top. 
	\end{multline*}
	Then observe that since we showed above that \(g_0(-, r)\) is disjoint from \(a_r\), we have 
	\begin{multline*}
		\bigvee_r g_0(-, r)^*
		\geq \bigvee_r a_r
		= \bigvee_r \bigvee_l \sbv{lr}{\delta > 0}\big(g_{m_l}(r + \delta, -) \wedge f(-, l\delta)\big)
		= \bigvee_l \sbv{lr}{\delta > 0} \bigvee_r \big(g_{m_l}(r + \delta, -) \wedge f(-, l\delta)\big)\\
		= \bigvee_l \sbv{lr}{\delta > 0}\left(f(-, l\delta) \wedge \bigvee_r \big(g_{m_l}(r + \delta, -) \big)\right)
		= \sbv{lr}{\delta > 0} \bigvee_l f(-, l\delta)
		= \top \qedhere
	\end{multline*}
\end{proof}

\begin{corollary}
	A vector lattice \(G\) is ru-Cauchy complete if and only if it is iru-complete.
\end{corollary}

\begin{proof}
	Suppose \(G\) is an ru-complete sub-vector lattice of \(H\). If \(g_n \to h\,(f)\) for some sequence \(\{g_n\} \subseteq G\) and regulator \(f \in G\) then \(\{g_n\}\) is ru-Cauchy and thus must converge to an element \(g_0 \in G\). Since ru-limits are unique, \(h = g_0 \in G\). We have shown \(G\) to be iru-complete.
	
	On the other hand, suppose \(G\) is iru-complete and consider a sequence \(\{g_n\} \subseteq G\) which is ru-Cauchy with regulator \(f \in G\). Identify \(G\) with the sub-vector lattice of \(\mcal{R}M\) which is  part of the (pointfree) Yosida representation (Theorem \ref{Thm:3}). By Theorem \ref{Thm:12} there is some \(h \in \mcal{R}M\) such that \(g_n \to h\,(f)\), and since \(G\) is iru-complete, \(h \in G\). We have shown \(G\) to be ru-Cauchy complete. 
\end{proof}

\begin{corollary}
	\(\mcal{R}M\) is iru-complete.
\end{corollary}

\subsection*{The iru-completion of a vector lattice}

\begin{definition*}[\(\mcal{T}G\), \(\tau_G\), \(\rho\)]
	We denote the insertion \(G \to \iruc_{\mcal{R}M} G\) by \(\map{\tau_G}{G}{\mcal{T}G}\), and we denote the insertion  \(\mcal{T}G \to \mcal{R}M\) by \(\rho\).
	\[
		\begin{tikzcd}
			G				\arrow{rr}{\mu_G}
							\arrow{dr}[swap]{\tau_G}
			&&\mcal{R}M\\
			&\mcal{T}G		\arrow{ur}[swap]{\rho}
		\end{tikzcd}
	\]
\end{definition*}

\begin{theorem}[\cite{BallHager:1999b}]\label{Thm:13}
	\(\mbf{irucV}\) is a full monoreflective subcategory of \(\mbf{V}\), and a reflector for the vector lattice \(G\) is the map \(\map{\tau_G}{G}{\mcal{T}G}\). 
\end{theorem}

\begin{proof}
	First note that \(\mcal{T}G\) is ru-complete because \(\mcal{R}M\) is ru-complete by Theorem \ref{Thm:12}. We appeal to the basic adjunction of Theorem \ref{Thm:3}.  Given a test \(\mbf{V}\) homomorphism \(\map{\theta}{G}{H}\) with \(H\) relatively uniformly complete, let \(m \equiv \mcal{K} \theta\) be the pointed frame homomorphism which realizes \(\theta\).
	\[
		\begin{tikzcd}
			G				\arrow{rr}{\mu_G}
							\arrow{dd}[swap]{\theta}
							\arrow{dr}[swap]{\tau_G}
			&&\mcal{R} M
							\arrow{dd}{\mcal{R}m}
			&M				\arrow{dd}{m}\\
			&\mcal{T}G
							\arrow{ur}[swap]{\rho}
			&&\\
		H					\arrow{rr}[swap]{\mu_H}
			&&\mcal{R}\mcal{K}H
			&\mcal{K} H
		\end{tikzcd}
	\]
	Since \(\mcal{R} m\) is relatively uniformly continuous by Proposition \ref{Prop:22}, it maps \(\mcal{T}G\) into \(\iruc_{\mcal{R K} H} \widehat{H}\), and since \(H\) is ru-complete,  \(\iruc_{\mcal{RK} H} \widehat{H} = \widehat{H}\). The desired factorization of \(\theta\) is \(\mu_H^{-1} \circ \mcal{R} m \circ \rho \circ \tau_G\). The factor \(\mu_H^{-1} \circ \mcal{R} m \circ \rho\) is unique because \(\tau\) is epic.
\end{proof}

\subsection*{\(G^* = \mcal{R}^*M\) when \(G\) is ru-Cauchy complete}

Throughout this subsection we assume \(G\) to be a vector lattice with standard compactification \(\map{q}{K}{M}\), identified with its image in the standard representation as a vector lattice in \(\mcal{E}q\) (Theorem \ref{Thm:3}). 
	
\begin{lemma}\label{Lem:10}
	Let \(\map{o_d}{K}{O_d}\) be the open quotient of the domain of reality \(d = f(-, -)\) of an element \(f \in G^+\). Then \(q\) factors through \(o_d\), say \(q = m_d \circ o_d\).   	    
	\[
	\begin{tikzcd}
		\mcal{O}\exR		\arrow{rr}{f'}
									\arrow{dd}[swap]{p}
		&&K							\arrow{dd}{q}
									\arrow{dl}[swap]{o_d}\\
		&O_d						\arrow{dr}
									\arrow{dr}[swap]{m_d}
		&\\
		\mcal{O}\mbb{R}				\arrow{rr}[swap]{f}
									\arrow{ur}{f''}
		&&M
	\end{tikzcd}
	\]
	Furthermore \(f\) factors through \(m_d\), say \(f = m_d \circ f''\), and each bounded element \(h \in \mcal{E}q\) factors through \(m_d\), say \(h = m_d \circ h''\). 
\end{lemma}	

\begin{proof}
	Because \(G\) is snugly embedded in \(\mcal{E}q\), \(q\) is the pushout of the open quotients arising from the domains of reality of the members of \(G\) and therefore factors through them all. 
	
	We define \(f''(U) \equiv o_d \circ f'(V)\) for any \(V \in \mcal{O}\ol{\mbb{R}}\) such that \(p(V) = U\). The function \(f''\) is well defined because \(p\) is surjective, and because for any \(V_i \in \mcal{O}\ol{\mbb{R}}\) we have
	\begin{multline*}
		p(V_1) = p(V_2)
		\implies\\
		 m_d \circ o_d \circ f'(V_1)
		= q \circ f'(V_1)
		= f \circ p(V_1)
		= f \circ p(V_2)
		= q \circ f'(V_2)
		= m_d \circ o_d \circ f'(V_2) \implies\\
		o_d \circ f'(V_1) 
		= o_d \circ f'(V_2).
	\end{multline*}
	The second implication holds because \(m_d\) is monic by virtue of being dense, and it is dense because \(q\) is dense and \(o_d\) is onto. Lastly observe that \(m_d \circ f''(U) = m_d \circ o_d \circ f'(V) = q \circ f'(V) = f \circ p(V) = f(U)\) for all \(U \in \mcal{O}\mbb{R}\).
	
	Since each bounded element \(h \in \mcal{E}^*q\) has an extension \(h'\) satisfying \(h \circ p = q \circ h'\), each truncation \(h \wedge n\) factors through \(q\) by Lemma \ref{Lem:21}(2), but since \(h\) is bounded, one of those truncation is \(h\) itself, which is only to say that \(h\) factors through \(q\), say \(h = q \circ h'\). Then \(h'' \equiv o_d \circ h'\) satisfies \(m_d \circ h'' = m_d \circ o_d \circ h' = q \circ h' = h\).
%
\end{proof}

\begin{lemma}\label{Lem:20}
	Let \(\map{c_b}{N}{\upset{b}} \equiv C_b = (a \mapsto b \vee a)\) be the closed quotient map of an element \(b\) of a frame \(N\). Then for any functions \(g,h \in \mcal{R}N\), 
	\[
		c_b \circ h 
		= c_b \circ g
		\implies \coz \abs{g - h}
		= \abs{g - h}(0, -) 
		\leq b.
	\]
\end{lemma}

\begin{proof}
	In view of the fact that 
	\[
		\abs{g - h}(0, -)
		= \sbv{lr}{\abs{U - V} \subseteq (0, -)} \big(g(U) \wedge h(V)\big)
		= \sbv{lr}{U \cap V = \emptyset}\big(g(U) \wedge h(V)\big),
	\]
	we get that 
	\begin{multline*}
		b \vee \abs{g - h}(0, -)
		= b \vee \sbv{lr}{U \cap V = \emptyset}\big(g(U) \wedge h(V)\big)
		= \sbv{lr}{U \cap V = \emptyset}\big(b \vee (g(U) \wedge h(V))\big)
		= \sbv{lr}{U \cap V = \emptyset}\big((b \vee g(U)) \wedge (b \vee h(V))\big)\\
		= \sbv{lr}{U \cap V = \emptyset}\big((b \vee h(U)) \wedge (b \vee h(V)))\big)
		= \sbv{lr}{U \cap V = \emptyset}\big(b \vee (h(U) \wedge h(V))\big)
		= \sbv{lr}{U \cap V = \emptyset}b 
		= b \qedhere 
	\end{multline*}
\end{proof}

\begin{proposition}\label{Prop:25}
	Assume the terminology of Lemma \ref{Lem:10}. If \(G\) is ru-Cauchy complete then for each bounded element \(h \in \mcal{R}^* O_d\) there is an element \(g_0 \in G^*\) such that \(g_0'' = h\).   
\end{proposition}

\begin{proof}
	By replacing \(f\) with \(f \vee 1\) if necessary, a replacement which does not alter its domain of reality, we may assume that \(f \geq 1\). Let \(f''\) be the element of \(\mcal{R}O_d\) mentioned in Lemma \ref{Lem:10} such that \(f'' \circ p = o_d \circ f'\). 	 
	\[
		\begin{tikzcd}
			\mcal{O}\exR					\arrow{r}{f'}
			\arrow{d}[swap]{p}
			&K								\arrow{d}{o_d}\\
			\mcal{O}\mbb{R}					\arrow{r}[swap]{f''}
			&O_d
		\end{tikzcd}	
		\qquad \qquad
%
	\begin{tikzcd}
		\mcal{O}\ol{\mbb{R}}		\arrow{rr}{g'}
		\arrow{dd}[swap]{p}
		&&K							\arrow{dd}{q}
		\arrow{dl}{o_d}\\
		&O_d 						\arrow{dr}{m_d}
		&\\
		\mcal{O}\mbb{R}				\arrow{rr}[swap]{g}
		\arrow[bend left]{uurr}{g'}
		\arrow{ur}[swap]{g''}
		&&M
	\end{tikzcd}
	\]
	Note the following.
	\begin{itemize}
		\item 
		Each bounded element \(g \in G^*\) factors through \(q\) by Lemma \ref{Lem:21}(2), say \(g = q \circ g'\).
		
		\item 
		\({G^*}' \equiv \setof{g'}{g \in G^*} = \mcal{R}K\) as a consequence of the Stone Weierstrass Theorem. The theorem is applicable because \(G^*\) is uniformly complete and its localic counterpart separates the points of its Yosida space, whose topology is (isomorphic to) \(K\).
		
		\item 
		For each \(g \in G^*\), \(o_d \circ g'\) is the element \(g''\) given in Lemma \ref{Lem:10} which satisfies \(m_d \circ g'' = g\).  If we let \({G^*}'' \equiv \setof{g''}{g \in G^*}\), we aim to show that \({G^*}'' = \mcal{R}^* O_d\).  
	\end{itemize} 	
	
	Now add a test frame map \(h\) to the diagram, assuming that \(h\) is bounded, say \(\abs{h} < t\) for \(0 < t \in \mbb{R}\). For each integer \(n\) let \(b_n \equiv f'(n, -]\) and let \(\map{c_n}{K}{\upset{b_n}} \equiv C_n = (a \mapsto a \vee b_n)\) be the associated closed quotient map.
	\[
		\begin{tikzcd}
			\mcal{O}\exR				\arrow{r}{f'}
										\arrow{d}[swap]{p}
			&K							\arrow{dr}{c_n}
										\arrow{d}{o_d}
			&\\
			\mcal{O}\mbb{R}				\arrow{r}[swap]{h}
										\arrow{ur}{g_n'}
			&O_d						\arrow{r}[swap]{k_n}
			&C_n
		\end{tikzcd}
	\]
	We claim that \(c_n\) factors through \(o_d\). For to say that \(o_d(a_1) = o_d(a_2)\) is to say that \(d \to a_1 = d \to a_2\), which implies \(d \wedge a_1 = d \wedge a_2\), and we claim that this implies that \(c_n(a_1) = c_n(a_2)\). That is because \(d \vee b_n = \top\), hence
	\begin{align*}
		a_1 \vee b_n
		&= (a_1 \vee b_n) \wedge \top
		= (a_1 \vee b_n) \wedge (d \vee b_n)
		= (a_1 \wedge d) \vee b_n\\
		&= (a_2 \wedge d) \vee b_n
		= (a_2 \vee b_n) \wedge (d \vee b_n)
		= (a_2 \vee b_n) \wedge \top
		= a_2 \vee b_n.
	\end{align*} 
	Let \(\map{k_n}{O_d}{C_d}\) be the map satisfying \(k_n \circ o_d = c_n\). 
	
	By virtue of being a closed quotient of a compact frame, \(c_n\) is a \(C^*\)-quotient and so \(k_n \circ h\) factors through it, say \(k_n \circ h = c_n \circ g_n'\). The initial factor \(g_n'\) lies in \(\mcal{R} K\), which coincides with \({G^*}'\) as noted above. By replacing \(g_n'\) with \(((-t) \vee g_n) \wedge t\) if necessary, we may assume that \(\abs{g_n'} \leq t\). 
	
	We claim that \(g_n'' \to h\, (f'')\) in \(\mcal{R}O_d\).  For a given positive integer \(k\) let \(m_k\) be the least integer such that \(m_k \geq 2tk\); we shall show that \(\abs{g_n'' - h} \leq f''/k\) for all \(n \geq m_k\) by showing that
	\[
		\abs{g_n'' - h}(r,-) 
		\leq (f''/k)(r, -)
		= f''(kr, -) 
	\]
	for all \(r \in \mbb{R}\). For that purpose fix an integer \(n \geq m_k \geq 2tk\). The displayed condition clearly holds if \(r \geq 2t\), for \(\abs{g_n'' - h}(2t, -) = \bot\) since \(\abs{g_n'' - h} \leq 2t\). The condition also holds for \(r < 1/k\), for in that case \(f''(kr, -) = \top\) since \(f'' \geq 1\). If \(1/k \leq r < 2t\) then 
	\[
		\abs{g_n'' - h}(r, -)
		\leq \abs{g_n'' - h}(0, -)
		\leq f''(n, -)
		\leq f''(m_k, -)
		\leq f''(2tk, -)
		\leq f''(rk,-).
	\] 
	(The second inequality is an application of Lemma \ref{Lem:20}.) This proves the claim, from which it follows that \(\{g_n''\}\) is an ru-Cauchy sequence in \({G^*}''\) with regulator \(f''\), and therefore \(\{g_n\}\) is an ru-Cauchy sequence in \(G^*\) with regulator \(f\). Since \(G\) is ru-Cauchy complete, it must then contain a unique element \(g_0\) such that \(g_n \to g_0\,(f)\). Since \(\abs{g_n} \leq t\) for all \(n\) it follows that \(g_0\) is bounded and so factors through \(m_d\), say \(g_0 = m_d \circ g_0''\). Thus we have in  \(\mcal{R}O_d\) that \(g_n'' \to h\,(f'')\) and \(g_n'' \to g_0''\,(f'')\), from which we conclude that \(g_0'' = h\) by virtue of the uniqueness of regular uniform limits. 
\end{proof}

\begin{corollary}[{\cite[12.3]{BallHager:1999a}}]
	If \(G\) is ru-Cauchy complete then the open quotient of the domain of reality of each positive member of \(G\) is a \(C^*\)-quotient. 
\end{corollary}

\begin{proof}
	In terms of the proof of Proposition \ref{Prop:25}, the desired extension of \(h\) is \(g_0' \circ p_*\), where \(p_*\) designates the adjoint map of \(p\), given by the formula  \(\map{p_*}{\mcal{O}\mbb{R}}{\mcal{O}\exR}
	= \big(a \mapsto \bigvee_{p(b) \leq a} b\ \big)\). 
\end{proof}

\begin{theorem}[{\cite[12.3]{BallHager:1999a}}]\label{Thm:1}
	If \(G\) is ru-Cauchy complete then \(G^* = \mcal{R}^*M\).
\end{theorem}

\begin{proof}
	Let \(\ssetof{d_i}{I}\) be the family of domains of reality of positive elements of \(G\). For each \(i \in I\) let \(\map{o_i}{K}{O_i}\) be the open quotient of \(d_i\), and since \(q\) factors through \(o_i\), let \(m_i\) be the map such that \(m_i \circ o_i = q\).  Now the fact that \(q\) is the cointersection of the family \(\{o_i\}_I\) of quotients means that the sink \(\{O_i\xra{m_i} M\}_I\) is the pushout of the source \(\{K\xra{o_i}O_i\}_I\) in \(\mbf{Frm}\), and gives \(q\) as \(m_i \circ o_i\) for (any) \(i \in I\). It follows that the pushout  \(\{\mcal{R}O_i\xra{\mcal{R}m_i}\mcal{R}M\}_I\) of the \(\mbf{V}\)-source \(\{\mcal{R}K\xra{\mcal{R} o_i}\mcal{R}O_i\}_I\), namely the \(\mbf{V}\)-subobject \(H \subseteq \mcal{R}M\) generated by \(\bigcup_I \mcal{R}p_i(\mcal{R}O_i)\), has Madden frame (isomorphic to) \(M\). Now \(\{d_i\}_I\) is closed under binary meets, for if \(d_i = g_i(-,-)\) for \(g_i \in G^+\) then \(d_1 \wedge d_2 = (g_1 \vee g_2)(-,-)\). Since each member of \(H\) is generated by finitely many elements of \(\bigcup_I \mcal{R}p_i(\mcal{R}O_i)\), it follows that the bounded elements of \(H\) are generated by the bounded elements of \(\bigcup_I \mcal{R}p_i(\mcal{R}O_i)\), all of which belong to \(G\) by Proposition \ref{Prop:25}. We have shown that \(G^* = H^* \subseteq \mcal{R}^*M\), and since \(H\) is order dense in \(\mcal{R}M\) and \(G\) is ru-Cauchy complete, \(G^* = \mcal{R}^*M\) by the Stone Weierstrass Theorem. 
\end{proof}

\section{The oru-completion of a vector lattice}\label{Sec:oruComp}

In this section we continue to use the symbol $G$ to represent a vector lattice with standard compactification \(\map{q}{K}{M}\) and Yosida representation \(\map{\mu_G}{G}{\widehat{G}} \subseteq \mcal{E}q\)

\begin{proposition}\label{Prop:2}
	An epicomplete vector lattice is oru-complete. 
\end{proposition}

\begin{proof}
	This follows from Proposition \ref{Prop:1}.
\end{proof}

We shall prove the converse of Proposition \ref{Prop:2}.

\begin{lemma}\label{Lem:2}
	If \(G\) is oru-complete then \(\mu_G\) is surjective, i.e., \(\widehat{G} = \mcal{R}M\). 
\end{lemma}

\begin{proof}
	Because \(G\) is iru-complete, \(\widehat{G}^* = \mcal{R}^* M\) by Theorem \ref{Thm:1}. And in \(\mcal{R}M\), the inequality \(g - g \wedge n \leq g^2/k\) holds for all integers \(n \geq k/4\) because the inequality \(r - r \wedge n \leq r^2/4\) holds in \(\mbb{R}\) for all integers \(n \geq k/4\) (\cite[3.1.2]{BallWalters:2002}). That is to say that \((g \wedge n) \to g\,(g^2)\) for all \(g \in \mcal{R}^+ M\). 
\end{proof}

Lemma \ref{Lem:3} recalls the basic facts regarding division in \(\mcal{R}M\).

\begin{lemma}[{\cite[3.3]{BallWalters:2002}}]\label{Lem:3}
	An element \(g \in \mcal{R}^+M\) has a multiplicative inverse \(1/g\) if and only if \(\coz g = \top\). When it exists, the inverse is given by the formula
	\[
		(1/g)(U) 
		= g(\tilde{U}),\ U \in \mcal{O}\mbb{R},
		\qtq{where} \tilde{U} 	
		\equiv \setof{x \in \mbb{R}}{1/x \in U}.
	\]
\end{lemma}

\begin{lemma}\label{Lem:6}
	Suppose \(g\) is an element of \(\mcal{R}^+ M\) such that \(\coz g = \top\). Then \((ng \wedge 1) \to 1\,(1/g)\).
\end{lemma}

\begin{proof}
	Given a positive integer \(k\), we have for any integer \(n > k/4\) that \(1 - ng \wedge 1 \leq 1/kg\) holds for any \(g \in \mcal{R}^+M\) because \(1 - nr \wedge 1 \leq 1/kr\) holds for any \(r \in \mbb{R}^+\) (\cite[3.1.2]{BallWalters:2002}).
\end{proof}

\begin{proposition}\label{Prop:4}
	If \(G\) is oru-complete then its Madden frame \(M\) is a \(P\)-frame, i.e., each cozero element of \(M\) is complemented.
\end{proposition}

\begin{proof}
	Fix a cozero element \(a \in M\), and let \(c_a\) and \(o_a\) be the closed and open quotient maps associated with \(a\). Lemma \ref{Lem:2} permits the identification of \(G\) with \(\mcal{R}M\), so we can find an element \(g_a \in G^+\) such that \(a = \coz g_a = g_a(0, -)\). Let \(H \equiv \mcal{R}\upset{a} \times \mcal{R}\downset{a}\) with projections \(\pi_{\upset{a}}\) and \(\pi_{\downset{a}}\), and let \(\theta\) be the product map making the diagram commute.
	\[
		\begin{tikzcd}
			\mcal{R}M		\arrow{r}{\mcal{R}o_a}
							\arrow{d}[swap]{\mcal{R}c_a}
							\arrow{dr}{\theta}
			&\mcal{R}\downset{a}\\
			\mcal{R}\upset{a}		
			&\mcal{R}\upset{a} \times\mcal{R}\downset{a} \equiv H
							\arrow{u}[swap]{\pi_{\downset{a}}}
							\arrow{l}{\pi_{\upset{a}}}
		\end{tikzcd}		
	\]
	Since the functor \(\mcal{R}\) preserves limits, this diagram is the result of applying \(\mcal{R}\) to the diagram of Lemma \ref{Lem:5}; in particular, \(H\) is isomorphic to \(\mcal{R}M_a = \mcal{R}(\upset{a} \times \downset{a})\) and \(\theta\) is isomorphic to \(\mcal{R}x_a\). The formula for \(\theta\) is 
	\[
		\theta(g)
		= ((\mcal{R}c_a)(g), (\mcal{R}o_a)(g))
		=(c_a \circ g, o_a \circ g),
		\qquad g \in G.
	\]
	\(\theta\) is injective because \(x_a\) is injective, and \(x_a\) is injective because the congruences associated with \(o_a\) and \(c_a\) are disjoint. 	Now
	\begin{align*}
		\coz(\mcal{R}c_a)(g_a)
		&= \coz (c_a \circ g_a)
		= c_a \circ g_a(0,-)
		= c_a(a)
		= a
		= \bot_{\upset{a}}
		\qquad \text{and}\\
		\coz (\mcal{R}o_a)(g_a)
		&= \coz (o_a \circ g_a)
		= o_a \circ g_a(0,-)
		= o_a(a)
		= a
		= \top_{\downset{a}}
	\end{align*}
	From the first line we deduce that \(c_a \circ g_a\) is the \(0\) element of \(\mcal{R}\upset{a}\). From the second we infer from Lemma \ref{Lem:6} that \(n (o_a \circ g_a) \wedge 1 \to 1\,(k)\) in \(\mcal{R}\downset{a}\), where \(k\) is the multiplicative inverse of \(o_a \circ g_a\) in \(\mcal{R}\downset{a}\). Together we conclude that \(n\theta(g_a) \wedge 1 \to (0,1)\,((0,k))\) in \(H\). 
	
	Since \(G\) is oru-complete there must be an element \(g_a' \in G\) such that \(\theta(g_a') = (0, 1)\). Now 
	\begin{align*}
		\theta(g_a' \wedge (1 - g_a'))
		&= \theta(g_a') \wedge (\theta(1) - \theta(g_a'))
		= (0, 1) \wedge ((1, 1) - (0, 1))
		= (0, 0)
		= 0_H \ \text{and}\\
		\theta(g_a' \vee (1 - g_a'))
		&= \theta(g_a') \vee (\theta(1) - \theta(g_a'))
		= (0, 1) \vee ((1, 1) - (0, 1))
		= (1, 1)
		= 1_H, 
	\end{align*}
	and \(\theta\) is one-one, hence \(g_a' \wedge (1 - g_a') = 0\) and \(g_a' \vee (1 - g_a') = 0\), from which it follows that \(\coz g_a'\) and \(\coz(1 - g_a')\) are complementary elements of \(M\). Since  \(n\theta(g_a) \wedge 1 \xra{ru} (0, 1) = \theta(g_a')\) it follows that \(ng_a \xra{ru}g_a'\) by Proposition \ref{Prop:22}, and since \(\{ng_a \wedge 1\}\) is increasing, Corollary \ref{Cor:1} yields \(\coz g_a' = \bigvee_n \coz (ng_a \wedge 1) = \coz g_a = a\). We have shown that \(a\) is complemented in \(M\).
\end{proof}

We summarize.

\begin{theorem}[\cite{BallHager:1993}]\label{Thm:2}
	The following are equivalent for a vector lattice \(G\) with Madden frame \(M\).
	\begin{enumerate}
		\item 
		\(G\) is oru-complete.
		
		\item 
		\(M\) is a \(P\)-frame and \(\widehat{G} = \mcal{R}M\).
		
		\item 
		\(G\) is epicomplete.
	\end{enumerate}
	Thus the oru-complete objects constitute a full monoreflective subcategory of \(\mbf{V}\), and a reflector for \(G\) is the arrow \((\mcal{R}p) \circ \mu_G\), where \(\map{p}{M}{P}\) is the \(P\)-frame reflection of \(M\) (\cite{BallWWZenk:2010}).
\end{theorem}

\begin{proof}
	The equivalence of (1) and (2) is the content of Lemma \ref{Lem:2} and Proposition \ref{Prop:4}. The equivalence (2) and (3) is the content of Theorem 3.9 of \cite{Ball:2018}.
\end{proof}

\section{Concerning the adjunction of Theorem \ref{Thm:3}}\label{Sec:RemLit}

The classical Yosida space \(X\) of a vector lattice \(G\) is the set of values of the unit \(1\), i.e., the set of convex subspaces maximal with respect to omitting \(1\), made into a compact Hausdorff space by the hull kernel topology. In the seminal paper \cite{Yosida:1942}, Yosida represents a bounded vector lattice \(G\) as a subobject of \(\mcal{C}X\), the vector lattice of continuous real valued functions on \(X\), with the constant function \(1\) as designated unit. The space \(X\) appears in our development in the form of its topology \(K\), which is the domain of the standard compactification \(\map{q}{K}{M}\) of \(G\). 

The Madden frame of \(G\) is the frame \(M\) of its \(\mbf{V}\)-kernels. In their influential paper \cite{MaddenVermeer:1986}, Madden and Vermeer represent an arbitrary vector lattice \(G\) as a sub-vector lattice of \(\mcal{R}M\).  The frame \(M\) appears in our development as the codomain of the standard compactification of \(G\), and the compactification map \(q\) itself realizes the inclusion \(G^* \to G\). 

In their important article \cite{HagerRobertson:1977}, Hager and Robertson provide a representation of an arbitrary vector lattice \(G\) which is in some sense intermediate between these two. In it, \(G\) is represented as a vector lattice in \(\mcal{D}X\), where \(X\) is the Yosida space of \(G\) and \(\mcal{D}X\) is the set of continuous functions \(X \to \ol{\mbb{R}}\) with dense domains of reality. If, in our development, we denote \(\mu_G(g)\) by \((\hat{g}', \hat{g})\) then \(\hat{g}'\) is the frame counterpart of the Hager-Robertson representation of \(g\) as an element of \(\mcal{D}X\), and \(\hat{g}\) is the Madden-Vermeer representation of \(g\) as an element of \(\mcal{R}M\). Thus we see that \(\mcal{E}q\) corresponds to a (possibly proper) subset of \(\mcal{D}X\) (cf.\ Example \ref{Ex:1}).

\end{document}